\documentclass[a4paper,english]{article}
\title{ A note related to the CS decomposition and the BK inequality for discrete determinantal processes. } 
\author{Andr\'e Goldman}

\usepackage[english]{babel}  
\usepackage{indentfirst}  
                                               
\usepackage[latin1]{inputenc}  
\usepackage[T1]{fontenc}  
\usepackage{amssymb,amsmath,amscd}

\newtheorem{conjecture}{Conjecture}

\newtheorem{lemma}{Lemma}

\newtheorem{prop}[lemma]{Proposition}
\newtheorem{theorem}{Theorem}
\newtheorem{remark}{Remark}

\newenvironment{proof}[1][]%
{\noindent{\bf Proof{}#1:}\par\nobreak}{\nopagebreak\hspace*{\fill} $\square$\vspace*{1.4ex}\par}

\begin{document}

\maketitle

\begin{abstract} We prove that for a discrete determinantal process the
BK inequality occurs for increasing events generated by simple points.
We give also some elementary but nonetheless  appealing relationship between a discrete determinantal process and the well-known CS decomposition.

\end{abstract}
\section{Introduction.}
There exists an extensive  mathematical literature, in several theoretical and applied areas, related to determinantal point processes. Let us cite,  to mention a few recent applied works, 
\cite{br},  \cite{KT}, \cite{KV}, \cite{d}, \cite{mr}, \cite{cn}.  A good overview of the main conceptual basis and properties can be found in \cite{lyon}
 and in the bibliography therein.\\
 From the theoretical point of view,  determinantal point processes could be defined (in Bourbaki-like-spirit) in the
general  locally compact Polish spaces setting, as point processes associated to some locally square integrable, Hermitian, positive semidefinite, locally trace-class operators and thereafter specialize them for particular cases namely to discrete determinantal processes. Regarding the latter, the approach of \cite{lyon} which consists to construct such processes, first in the most elementary discrete context and then  gradually extend them to the general situation, provides, in our opinion, many advantages. It turns out also that some results for the most general processes are proved only \cite{golP}, \cite{lyon}, or more simply \cite{lyons}, indirectly from the corresponding results of the basic processes.\\
The basic elementary determinantal point process can be described via the exterior product concept, as follows.\\
Fix $1< p < N$ and let 
  $ \mathfrak{Z}=\{z^{1},...,z^{p}\} $, $ 1<p<N $, be a set of orthonormal vectors
  in $ \mathbb{C}^{N} $. Denote 
 
$$
z^{i}=(z^{i}_{1},\dots,z^{i}_{N})^{t}, i=1,\dots,p
$$
and
$$
z_{i}=(z^{1}_{i},\dots,z^{p}_{i}), i=1,\dots,N.
$$

The associated determinantal process $ \phi(\mathfrak{Z}) $
is a point process, view as a random subset of $\mathcal{N}=\{1,\dots,N\}$
of cardinality $ \vert \phi(\mathfrak{Z})\vert = p $, characterized \cite{lyo}, \cite{lyon}
by the formulas
\begin{equation}\label{det} 
P\{\{i_{1},\dots,i_{p}\} = \phi\}=
 \mid(\bigwedge_{i=1}^{p}z^{i})_{\{i_{1},\dots,i_{p}\}}\mid^{2}
= [\det((z_{i_{j}}^{k})_{k,j=1,\dots,p})]^{2}
\end{equation}
for all subsets $\{i_{1},\dots,i_{p}\} \subset \mathcal{N}$.
Note also that (\ref{det}) implies 
\begin{equation}
P\{\{i_{1},\dots,i_{k}\} \subset \phi\}=
\parallel \bigwedge_{j=1}^{k}z_{i_{j}}\parallel ^{2}
\end{equation}
for all $1\leq k \leq p$.

Let $ E=E(\mathfrak{Z})\subset
 \mathbb{C}^{N}$ 
 be 
the vector space spanned by  $ \mathfrak{Z} $. For all sets of linearly independent vectors $v^{i}\in E$, $i=1,\dots,p$,   we have $\bigwedge_{i=1}^{p}v^{i}=a\bigwedge_{i=1}^{p}z^{i}$ with $a\neq 0$ thus, in particular, if 
$ \tilde{\mathfrak{Z}}=\{\tilde{z}^{1},...,\tilde{z}^{p}\} $
is another orthonormal basis of $ E=E(\mathfrak{Z}) $ then  
$$
\mid(\bigwedge_{i=1}^{p}z^{i})_{\{i_{1},\dots,i_{p}\}}\mid
=
\mid(\bigwedge_{i=1}^{p}\tilde{z}^{i})_{\{i_{1},\dots,i_{p}\}}\mid
$$ 
for every
$\{i_{1},\dots,i_{p}\}\subset \mathcal{N}$ and consequently
$ \phi(\mathfrak{Z})=\phi(\tilde{\mathfrak{Z}}) $. \\
Remark also that if $ \mathfrak{Z}^{\perp}=\{z^{p+1},...,z^{N}\} $ is
an orthonormal basis of  the orthogonal complement 
$ E(\mathfrak{Z})^{\perp} $ of  $E(\mathfrak{Z})$ in $\mathbb{C}^{N}$
then obviously
$$
\phi(\mathfrak{Z}^{\perp})=\{1,\dots,N\}\setminus \phi(\mathfrak{Z}).
$$
A remarquable example of a non-trivial basic determinantal process is given by uniform spanning tree measure on a finite connected graph G. Roughly speaking, if G fixed, is arbitrary edge-oriented and M is the vertex-edge incidence matrix (the columns been indexed by vertexes) then the determinantal process associated to the vector space spanned by column vectors but one, provides a uniform probability on spanning trees. This result due \cite{bp} is called the Transfer Current Theorem. For more details, with clever short  proofs, see \cite{lyon} 2.6. p.8. 
Some extensions of this result are given in \cite{BLPS} with a serie of open questions and conjectures. Among them, Conjecture 4.6. related to BK-type inequality after
J. van den Berg and  H.Kesten.\\
 Recall that an event $\mathfrak{A}\subset 2^{\mathcal{N}}$, 
 $\mathcal{N}=\{1,\dots,N\}$,
 is called increasing if whenever $A\in \mathfrak{A}$ and 
 $n \in \mathcal{N}$, we have also $A\cup \{n\}\in \mathfrak{A}$. For a pair $\mathfrak{A}$, $\mathfrak{B}$ $\subset 2^{\mathcal{N}}$ of increasing events the disjoint intersection  $\mathfrak{A}\circ\mathfrak{B}$ is then defined \cite{BK} by
\begin{equation}
\mathfrak{A}\circ\mathfrak{B} 
=\{ K \subset \mathcal{N}:\exists \quad L\in \mathfrak{A}, \quad M \in \mathfrak{B}, \quad L,M \neq \emptyset, \quad L \cap M = \emptyset,
K \supset L\cup M
 \}.
 \end{equation}
A point process $\psi$ on $\mathcal{N}$ is said to have the van den Berg - Kesten property (in short the BK property) if 
\begin{equation}\label{BK}
P\{ \psi \in \mathfrak{A}\circ\mathfrak{B}\}
\leq 
P\{\psi \in \mathfrak{A}\}
\times
P\{\psi \in \mathfrak{B}\}
\end{equation}
for every pair of increasing events. In \cite{BK} J. van den Berg and  H.Kesten proved that inequality (\ref{BK}) is satisfied when $\psi$ is related to a product probability on $2^{\mathcal{N}}$. In the basic determinantal process setting the Conjecture 4.6. which  states that the same is true for the spanning tree determinantal point processes is still unsolved. The question of whether general determinantal processes have the BK property was raised in \cite{lyo}.

The purpose of this  note is twofold. First,  we introduce a new method to investigate discrete determinantal processes using the CS decomposition (CSD) of partitioned unitary matrix which is a useful non-trivial tool in numerical linear algebra: a precise statement of CSD is given in section 2.   We show that the CSD  gives a pertinent description of conditioning and provides (at least in our opinion...) a suggestive perspective for future investigations: see for example the result given by the proposition \ref{prop2} below which seems us new and the results of \cite{gol}. Furthermore, this should be  an appropriate framework for computational needs. \\ 
Second, we study the BK inequality. We prove that the BK inequality (\ref{BK}) is satisfied for all discrete determinantal processes when the increasing events
$\mathfrak{A}$ and $\mathfrak{B}$ are  generated by simple points: Theorem \ref{t2} in section 4 and section 5. We conjecture also that:
\begin{conjecture}\label{Con} 
For all $n\geq 2$ we have
\begin{equation}\label{C}
\begin{split}
& P\lbrace  A \not \subset \phi \mid A_{i} \not \subset \phi, \forall i=1,\dots,n\rbrace \\
& \leq P\lbrace A \not \subset \phi \mid A_{i} \not \subset \phi, \forall i=1,\dots,n-1\rbrace
 \end{split}
\end{equation}
for every choice $A$, $A_{i}$, $i=1,\dots,n$, of  disjoint subsets of $\{1,\dots, N\}$ such that 
$P\lbrace  A_{i} \not \subset \phi, \forall i=1,\dots,n\rbrace
> 0 $.
\end{conjecture}
If the Conjecture \ref{Con} holds then it can be shown that the BK inequality (\ref{BK}) is satisfied for increasing events
$\mathfrak{A}$ and $\mathfrak{B}$ generated by disjoint sets: Theorem \ref{t2g} in section 3. When the sets above are  reduced to being  simple points then  the inequality (\ref{C}) is a well-known result. 
For general sets, note that $P\lbrace  A \not \subset \phi \mid A_{1} \not \subset \phi, \rbrace \leq P\lbrace A \not \subset \phi \rbrace$,
$A\cap A_{1}=\emptyset $  is the classical correlation inequality \cite{lyo} and that (\ref{C}) was obtained in \cite{gol} for $n=2,3 $ with  the precise values of conditional probabilities. 
\begin{remark}\label{R0}
Note that for the process $\psi$  related to a product probability on $2^{\mathcal{N}}$ the counting
random variables $\vert \psi \cap A_{i} \vert$ (the sets $A_{i}$, $i=1,\dots,n$, being  disjoint) are independent and thus
the inequality (\ref{C}) becomes trivial. However, the situation is less obvious if the process $\psi$ is conditioned to have exactly $k$ points  $1<k<N$. In the particular case
when the conditioned process $\psi_{k}$ assigns equal probability to all subsets $\{i_{1},\dots,i_{k}\} \subset \mathcal{N}$,
that is if 
 $P\lbrace \psi_{k} = \{i_{1},\dots,i_{k}\}\rbrace
= 1/ \binom{N}{k}$,
it was proved in \cite{BJ} that  $\psi_{k}$ has the BK property. As regards to the inequality (\ref{C}) we have, 
with the choice  $P\lbrace  i \not \in \psi_{k},  i_{j} \not \in \psi_{k}, \forall j=1,\dots,n\rbrace>0$, 
 $$ 
 P\lbrace  i \not \in \psi_{k} \mid i_{j} \not \in \psi_{k}, \forall j=1,\dots,n\rbrace \\
 = \frac{\binom{N - n - 1}{k}}{\binom{N - n }{k}}
$$
and, consequently, the inequality (\ref{C}) is  equivalent, for simple points, to the well-known   log-concave inequality
$\binom{N - n - 1}{N-k}\times\binom{N - n + 1}{N-k} \leq \binom{N - n}{N-k}^{2}$ and thus is fulfilled. Likewise,  for general sets, the correlation inequality 
$P\lbrace  A \not \subset \psi_{k} \mid A_{1} \not \subset \psi_{k} \rbrace \leq P\lbrace A \not \subset \psi_{k} \rbrace$,
$A\cap A_{1}=\emptyset $, $\vert A \vert = n$, $\vert A_{1}\vert = m$ with (the non-trivial case) $n + m \leq k$,  follows from BK property and
is equivalent to the log-concave inequality $\binom{N}{N-k}\times\binom{N - n - m}{N-k} \leq \binom{N-n}{N-k}\times \binom{N - m}{N-k} $. For $n\geq 2$ it is easy to see that the validity of the inequality (\ref{C}) depends on whether or not the functions of the form
$$u\rightarrow \sum_{i_{1}=0}^{n_{1}}...\sum_{i_{M}=0}^{n_{M}}
\binom{N - u - (i_{1}+...+i_{M})}{N-k-M} $$
are log-concave, a question  which does not seem to me to have been really investigated. Finally, the occurrence of log-concave criteria for negative dependence properties is not quite a surprise, see for example \cite{Pe} .
\end{remark}
\section{ The CS decomposition and the basic determinantal point process}
Following \cite{pw}   the general   CS decomposition (CSD) for a matrix $Q$ from the unitary group $U(N)$ specifies that for any $2\times 2$ partitionning
\begin{equation}
\begin{split}
& \qquad \quad c_{1}\qquad  c_{2}\\
& Q= \begin{bmatrix}
Q_{11}& Q_{12}\\
Q_{21}& Q_{22}\\
\end{bmatrix}
\begin{matrix}
r_{1}\\
r_{2}\\
\end{matrix}
\end{split}
, \qquad N=r_{1} + r_{2}=c_{1} + c_{2},
\end{equation} 
there exist unitary matrices $U_{1}$, $U_{2}$, $V_{1}$, $V_{2}$
such that (here all unnamed blocks of the matrices below  are always zero and the superscript $H$ represents the conjugate transpose)
\begin{equation}\label{CSD} 
\begin{bmatrix}
U_{1}^{H}& \\
 & V_{1}^{H}\\
\end{bmatrix}
Q
\begin{bmatrix}
U_{2} &\\
  & V_{2} \\
\end{bmatrix}
=
\left [\begin{array}{c|c}
U_{1}^{H}Q_{11}U_{2}& U_{1}^{H}Q_{12}V_{2} \\
\hline
V_{1}^{H}Q_{21}U_{2}& V_{1}^{H}Q_{22}V_{2} \\
\end{array}\right]
= \overset{c_{1} \quad c_{2}} 
{ \left[ \begin{array}{c|c}
D_{11}& D_{12} \\
\hline
D_{21}& D_{22} \\
\end{array}\right]}
\begin{matrix}
r_{1}\\
r_{2}\\
\end{matrix}
\end{equation}
where
  the matrices
$$
D_{11}=
\begin{bmatrix}
I & &\\
& C & \\
& & \mathbf{0}_{c}\\
\end{bmatrix},
\quad 
D_{12}=
\begin{bmatrix}
\mathbf{0}_{s}^{H} &  &\\
& S & \\
& & I\\
\end{bmatrix},
\quad 
D_{21}=
\begin{bmatrix}
\mathbf{0}_{s} &  &\\
& S & \\
& & I\\
\end{bmatrix},
$$
$
D_{22}=
\begin{bmatrix}
I &  &\\
& -C & \\
& & \mathbf{0}_{c}^{H}\\
\end{bmatrix}
$
are diagonal
with
$
C\equiv diag(\cos\theta_{1},\dots,\cos\theta_{s}),
$
\\
$S \equiv diag(\sin\theta_{1},\dots,\sin\theta_{s})$,
$1>\cos\theta_{1}>\dots > \cos\theta_{s}>0$.
\\ 
In some cases the matrices of zeros $\mathbf{0}_{s}$ and $\mathbf{0}_{c}$ as well as the unit matrices $I$ could be nonexistent. See \cite{pw} THEOREM 1 (page $12$ and the discussion that follows) for the full statement, and below for a detailled description given  from Jordan's geometrical point of view.\\
The CS decomposition is a deep result which has a long history going back to the work of Camille Jordan in 1875 on angles between subspaces in $\mathbb{R}^{n}$ \cite{jor}. Nowaday it is a popular tool in numerical linear algebra, useful for solving various questions as, for example, constrained least square problems, computing principal angles between subspaces, the generalized singular value decomposition, quantum computing and more else 
\cite{bai}, \cite{sut}, \cite{gvl}, \cite{cs},  \cite{pw}.\\
Now, let $ E\subset \mathbb{C}^{N}$ be a vector space of dimension $ 1 < p< N $, $ \mathfrak{Z}=\{z^{1},...,z^{p}\} $
an orthonormal basis of  $ E $ and
$ \mathfrak{Z}^{\perp}=\{z^{p+1},...,z^{N}\} $  an orthonormal basis of the orhogonal complement
$ E(\mathfrak{Z})^{\perp} $. Fix $1\leq n\leq p$ and consider the CSD of the partitioned  unitary matrix $Q=(z^{1},...,z^{p},
z^{p+1},...,z^{N}) $ 
\begin{equation}
\begin{split}
& \qquad \quad p \quad  N-p\\
& Q= \begin{bmatrix}
Q_{11}& Q_{12}\\
Q_{21}& Q_{22}\\
\end{bmatrix}
\begin{matrix}
n\\
N-n\\
\end{matrix}
\end{split}
\end{equation}
It follows from (\ref{CSD}) that the columns vectors of these two matrices
\begin{equation}\label{z} 
\begin{bmatrix}
U_{1}D_{11}\\
V_{1}D_{21}\\
\end{bmatrix}
, \quad
\begin{bmatrix}
U_{1}D_{12}\\
V_{1}D_{22}\\
\end{bmatrix}
\end{equation}
are respectively orthonormal basis of  $ E $ and 
$ E(\mathfrak{Z})^{\perp} $.\\
Now we will detail different  cases given by (\ref{z})  which need to  be distinguished. 
The description given here is somewhat lengthy but, in our opinion,
useful for both theoretical and computational purpose.
We note by 
$ e(k) $, $k=1,\dots,N $,  the nul vector of the space 
$\mathbb{C}^{k}$. Note also the slight change with regard to angles appearing in CSD (\ref{CSD}) which allows values $0$ and $\pi /2$ in order to recover all Jordan's principal angles.
\begin{enumerate}
\item[I.-] The case $ n<p $ and $p+n < N$.\\
There exist:
\begin{enumerate} 
\item
a sequence $ u^{1},\dots,u^{n} $ of orthonormal vectors in
$\mathbb{C}^{n}$,
\item three sequences of 
mutually orthonormal vectors in $\mathbb{R}^{N-n}$  \\ 
$ \mathfrak{V}=\{V^{1},\dots,V^{n}\} $,
$\mathfrak{W}= \{W^{1},\dots,W^{p-n}\} $ and \\
$ \mathfrak{\tilde{W}}=\{\tilde{W}^{1},\dots,\tilde{W}^{N-p-n}\} $,
\item Jordan angles $0\leq \theta_{1}\leq \dots  \leq \theta_{n}\leq \pi/2 $
\end{enumerate}
such that noting
\begin{enumerate}
\item[]
$
z^{i}=
\begin{bmatrix}
u^{i}cos\theta_{i}\\
V^{i}sin\theta_{i}
\end{bmatrix}
  i=1,\dots,n, \quad
z^{i}=
\begin{bmatrix}
e(n)\\
W^{i}
\end{bmatrix}
 \quad i=n+1,\dots,p,
$
  
\item[]
$
z^{p+i}=
\begin{bmatrix}
u^{i}sin\theta_{i}\\
-V^{i}cos\theta_{i})
\end{bmatrix}
  i=1,\dots,n,
$
\item[]
$
z^{p+n+i}=
\begin{bmatrix}
e(n)\\
\tilde{W}^{i}
\end{bmatrix}
i=1,\dots,N-p-n,
$
\end{enumerate}
the sequence $ \mathfrak{Z}=\{z^{1},...,z^{p}\} $ is an orthonormal basis of  $ E $  and the sequence $ \mathfrak{Z}=\{z^{p+1},...,z^{N}\} $ is an orthonormal basis of the orthogonal complement $ E^{\perp} $.  
\item[II.-]The case $ n < p $ and $ p+n > N $.\\ 
There exist:
\begin{enumerate}
\item
a sequence $ u^{1},\dots,u^{n} $ of orthogonal vectors in
$\mathbb{C}^{n}$,
\item two sequences of mutually orthogonal vectors in $\mathbb{C}^{N-n}$  \\
$ \mathfrak{V}=\{V^{1},\dots,V^{N-p}\} $ and
$ \mathfrak{W}=\{W^{1},\dots,W^{p-n}\} $,  
\item 
Jordan angles 
$0 =\theta_{1}=\dots=\theta_{n+p-N}
\leq \dots  \leq \theta_{n}\leq \pi/2 $
\end{enumerate}
such that noting
\begin{enumerate}
\item[]
$
z^{i}=
\begin{bmatrix}
u^{i}\\
e(N-n)
\end{bmatrix}
i=1,\dots,n+p-N,
$
 
\item[]
$
z^{i}=
\begin{bmatrix}
u^{i}cos\theta_{i}\\
V^{i-n-p +N}sin\theta_{i})
\end{bmatrix}
i=n+p-N +1,\dots,n,
$
\\
$
z^{n+i}=
\begin{bmatrix}
e(n)\\
W^{i}
\end{bmatrix}
i=1,\dots,p-n,
$
\item[]
$
z^{p+i}=
\begin{bmatrix}
u^{n+p-N+i}sin\theta_{n+p-N+i}\\
-V^{i}cos\theta_{n+p-N+i}
\end{bmatrix}
i=1,\dots,N-p,
$
\end{enumerate}
the set $ \mathfrak{Z}=\{z^{1},...,z^{p}\} $ is an orthonormal basis of $ E $  and the set $ \mathfrak{Z}=\{z^{p+1},...,z^{N}\} $ is an orthonormal basis of $ E^{\perp} $.  
\item[III.-]The case $ n < p $ and $ p+n = N. $\\ 
There exist:
\begin{enumerate}
\item
a sequence $ u^{1},\dots,u^{n} $ of orthogonal vectors in
$\mathbb{C}^{n}$,
\item two sequences of mutually orthogonal vectors in $\mathbb{C}^{N-n}$ \\
$ \mathfrak{V}=\{V^{1},\dots,V^{n}\} $ and
$ \mathfrak{W}=\{W^{1},\dots,W^{p-n}\} $,  
\item 
Jordan angles 
$0 \leq \theta_{1}\leq \dots \leq \theta_{n}\leq \pi/2 $
\end{enumerate}
such that noting
\begin{enumerate}
\item[]
$
z^{i}=
\begin{bmatrix}
u_{i}cos\theta_{i}\\
V^{i}sin\theta_{i}
\end{bmatrix}
i=1,\dots,n,
\quad
z^{n+i}=
\begin{bmatrix}
e(n)\\
W^{i}
\end{bmatrix}
i=1,\dots,p-n,
$

\item[]
$
z^{p+i}=
\begin{bmatrix}
u^{i}sin\theta_{i}\\
-V^{i}cos\theta_{i}
\end{bmatrix}
i=1,\dots,n,
$
\end{enumerate}
the set $ \mathfrak{Z}=\{z^{1},...,z^{p}\} $  is an orthonormal basis of de $ E $  
and the set $ \mathfrak{Z}=\{z^{p+1},...,z^{N}\} $ is an orthonormal basis of $ E^{\perp}. $  
\item[IV.-]The case $ n = p $.\\ 
With the notations of points I-III:
\begin{enumerate}
\item[(i)]$ 2p< N $
\begin{enumerate}
\item[]
$
z^{i}=
\begin{bmatrix}
u^{i}cos\theta_{i}\\
V^{i}sin\theta_{i}
\end{bmatrix}
i=1,\dots,p, \quad
z^{p+i}=
\begin{bmatrix}
u^{i}sin\theta_{i}\\
-V^{i}cos\theta_{i}
\end{bmatrix}
 i=1,\dots,p, 
 $
\item[]
$
z^{2p+i}=
\begin{bmatrix}
e(n)\\
W^{i}
\end{bmatrix}
i=1,\dots,N-2p.
$
\end{enumerate}
\item[(ii)]$ 2p > N $
\begin{enumerate}
\item[]
$
z^{i}=
\begin{bmatrix}
u^{i}\\
e(N-n)
\end{bmatrix}
i=1,\dots,2p-N,
$
\item[]
$
z^{i}=
\begin{bmatrix}
u^{i}cos\theta_{i}\\
V^{i-2p+N}sin\theta_{i}
\end{bmatrix}
i=2p-N+1,\dots,p,
$
\item[]
$
z^{i}=
\begin{bmatrix}
u^{i+p-N}sin\theta_{i+p-N}\\
-V^{i-p}cos\theta_{i+p-N}
\end{bmatrix}
i=p+1,\dots,N.
$
\end{enumerate}
\item[(iii)]$ 2p = N $
\begin{enumerate}
\item[]
$
z^{i}=
\begin{bmatrix}
u^{i}cos\theta_{i}\\
V^{i}sin\theta_{i},
\end{bmatrix}
 \quad
z^{p+i}=
\begin{bmatrix}
u^{i}sin\theta_{i}\\
-V^{i}cos\theta_{i}
\end{bmatrix}
 i=1,\dots,p. 
$
\end{enumerate}
\end{enumerate}
\end{enumerate}

By reordering the rows of $Q$  the  procedure described above works for every subset 
$J=\{x_{1},\dots,x_{n}\}\subset \{1,\dots,N\}  $, $ 1\leq n\leq p $,
and gives related basis of the spaces $E$ and $ E^{\perp} $.
Note that in the Euclidean context, that is for $E\subset \mathbb{R}^{N}$ and the CSD applied to orthogonal matrices,
the angles appearing in the  CS decomposition (related to $J$) are principal Jordan angles between
the space E and the basic subspace
$$
\mathbb{R}_{J}^{N}=\{x=(x_{k})\in\mathbb{R}^{N};x_{k}=0 \quad
if \quad k\notin J\}.
$$

An important statistical application of principal angles is the canonical correlation analysis (CCA) of H.Hotelling \cite{hot}. In order to develop a unified algebraic formulation of concepts in multivariate analysis (as for example CCA), S.Afriat   has thoroughly studied in \cite{afr} (see also \cite{m-BI}) the geometry of subspaces in $\mathbb{R} ^{N}$  in terms of orthogonal and oblique projectors and has introduced, among others, the notation of so-called \textit{multiplicative cosine and sine }     
$$ 
\cos\{E,\mathbb{R}_{J}^{N}\}=\prod_{i=1}^{n}\cos\theta_{i}, 
$$
$$ 
\sin\{E,\mathbb{R}_{J}^{N}\}=\prod_{i=1}^{n}\sin\theta_{i}. 
$$

The basis of $E$ given by the CSD is a pertinent tool for the study of
the associated determinantal process. For example it gives at once
\begin{prop}\label{prop1}
For a set  $ J=\{x_{1},\dots,x_{n}\} $, $ n \leq p $  we have:
\begin{enumerate}
 \item 
\begin{equation}\label{it}
P\{\vert J\cap \phi\vert = n\}=\prod_{i=1}^{n}\cos^{2}\theta_{i}
\end{equation}
and for
$k=1,\dots, n-1$
\begin{equation}\label{it1}
P\{\vert J\cap \phi\vert = k\}=\sum_{1\leq i_{1}<\dots < i_{k}\leq n}\prod_{j=1}^{k}\cos^{2}\theta_{i_{j}}\times
\prod_{j\not \in\{i_{1},\dots , i_{k}\} }\sin^{2}\theta_{j}.
\end{equation}
\item
\begin{equation}\label{itt}
P\{\vert J\cap \phi^{c }\vert = n\}=\prod_{i=1}^{n}\sin^{2}\theta_{i}
\end{equation}
and for $k=1,\dots, n-1 $
\begin{equation}\label{it2}
P\{\vert J\cap \phi^{c }\vert = k\}=\sum_{1\leq i_{1}<\dots < i_{k}\leq n}\prod_{j=1}^{k}\sin^{2}\theta_{i_{j}}\times
\prod_{j\not \in\{i_{1},\dots , i_{k}\} }\cos^{2}\theta_{j}.
\end{equation} 
\item  
If $n < p$ and  $ P\{ J \subset\phi\} >0$ then the conditioned process $\{\phi\vert \ J\subset\phi\}\setminus J$ is determinantal such that 
$\{\phi\vert \ J\subset\phi\}\setminus J=
\phi(\mathfrak{W}).$
\item
If $N-p > n$ and  $ P\{ J \subset\phi^{c}\} >0$ then the conditioned process $\{\phi\vert \ J\subset\phi^{c}\}$ is determinantal such that 
$\{\phi\vert \ J\subset\phi^{c}\}=
\phi(\mathfrak{V}\cup \mathfrak{W}).$
\item
If $ P\{ J \subset\phi\} >0$ 
then for all $ K \subset \{1,\dots,N\}\setminus J $ we have
\begin{equation}\label{it3}
P\{K\subset \phi(\mathfrak{W})\} 
\leq P\{K\subset \phi\},
\end{equation}
if $ P\{ J \subset\phi^{c}\} >0$ then
\begin{equation}\label{it4} 
P\{K\subset \phi\} 
\leq P\{K\subset \phi(\mathfrak{V}\cup \mathfrak{W})\}.  
\end{equation}
\end{enumerate}
\end{prop}
\begin{remark}\label{R1}
The fact that the conditioned processes 
$\{\phi\vert \ J\subset\phi\}\setminus J$ 
and $\{\phi\vert \ J\subset\phi^{c}\}$ are determinantal, as well
as inequalities (\ref{it3}) and (\ref{it4}), are well-known results proved by R.Lyons \cite{lyo}.
\end{remark}
\begin{remark}\label{R2}
Regarding items 1 and 2 of Proposition \ref{prop1}, it was proved more generally
in \cite{hkpv} Theorem 5, that for general determinantal processes with  trace-class (both discrete and continuous case)  kernels, the number of points
in the process has the distribution of a sum of independent Bernoulli random variables.
\end{remark}
A more elaborate information can be obtained from this point of view. 
For example:
\begin{prop}\label{prop2}
Consider the discrete determinantal process  $ \phi = \phi(\mathfrak{Z}) $ 
associated to a set $ \mathfrak{Z}=\{z^{1},...,z^{p}\} $, $ 1<p<N $, 
 of orthonormal vectors in $ \mathbb{C}^{N} $. Fix points
$J=\{x_{1},\dots,x_{n}\}\subset \{1,\dots,N\}  $, $ 1\leq n\leq p $,
such  that $P\{\{x_{2},\dots,x_{n}\}\subset \phi^{c}\} > 0 $.
With the choice (to simplify the notations) $x_{i}=i $, $ i=1,\dots,n $, we have
\begin{eqnarray}\label{pp}
\lefteqn{ \bigg\vert\langle z_{1},z_{n} \ \rangle +
\sum_{k=1}^{n-2}(-1)^{k}\sum_{2\leq i_{1}<\dots < i_{k}\leq n-1}
  \langle z_{1}\wedge(\bigwedge_{j=1}^{k} z_{i_{j})},
 z_{n}\wedge(\bigwedge_{j=1}^{k} z_{i_{j}}) \ \rangle
\bigg\vert^{2}} \nonumber \\
& & {}=
P\{\{x_{2},\dots, x_{n}\}\subset \phi^{c}\} \times 
P\{\{x_{2},\dots, x_{n-1}\}\subset \phi^{c}\}\\
& & {}\quad \times
\bigg[ P\{x_{1}\in \phi \vert \{x_{2},\dots,x_{n}\}\subset \phi^{c}\}
-P\{x_{1}\in \phi \vert \{x_{2},\dots,x_{n-1}\}\subset \phi^{c}\}\bigg]. 
\nonumber
\end{eqnarray}
\end{prop}
\begin{proof} The left and rigt sides of (\ref{pp})
do not depend of the choice of the basis of $E$. Choose the basis given by the CS decomposition related to the set $J=\{2,\dots,n-1\}$, with the reordering ($2,\dots,n-1,1,n)^{t}$
and $N-p - n +2>0$ (the general situation, case I). The  first n-coordinates of these  basis
 have  the following form:
\begin{equation}\label{pr1}
\begin{split}
\begin{matrix}
2\\
\vdots\\
n-1\\
1\\
n
\end{matrix} 
& \left[ \begin{matrix}
\cos \theta_{1}u^{1}_{1}  \dots  \cos \theta_{n-2}u^{n-2}_{1}
  &  0 \quad \dots \quad 0 \\
\vdots  \qquad & \vdots    \\
 \cos \theta_{1}u^{1}_{n-2}  \dots  \cos \theta_{n-2}u^{n-2}_{n-2}
  &   0 \quad \dots \quad 0 \\
\sin\theta_{1}V^{1}_{1}   \dots  \sin \theta_{n-2}V^{n-2}_{1} & W^{1}_{1} \dots W^{p-n+2}_{1}\\
\sin\theta_{1}V^{1}_{2}   \dots  \sin \theta_{n-2}V^{n-2}_{2} & W^{1}_{2} \dots W^{p-n+2}_{2}
\end{matrix} \right. 
\\
& \\
&
\qquad   
 \left. \begin{matrix}
\sin \theta_{1}u^{1}_{1}  \dots  \sin \theta_{n-2}u^{n-2}_{1}
  &  0 \quad \dots \quad 0 \\
\vdots  \qquad & \vdots    \\
 \sin \theta_{1}u^{1}_{n-2}  \dots  \sin \theta_{n-2}u^{n-2}_{n-2}
  &  0 \quad \dots \quad 0 \\
-\cos\theta_{1}V^{1}_{1}   \dots  -\cos \theta_{n-2}V^{n-2}_{1} & \tilde{W}^{1}_{1}  \dots \  \tilde{W}^{N+2-n-p}_{1}\\
-\cos\theta_{1}V^{1}_{2}   \dots  -\cos \theta_{n-2}V^{n-2}_{2} & \tilde{W}^{1}_{2}  \dots   \tilde{W}^{N+2-n-p}_{2}
\end{matrix} \right]   
\end{split} 
\end{equation} 
It follows from Proposition \ref{prop1} that
\begin{enumerate}
\item[(a)]
$P\{\{x_{2},\dots, x_{n-1}\}\subset \phi^{c}\}=\prod_{i=1}^{n-2}
\sin^{2}\theta_{i}$
\item[(b)]
$P\{\{x_{2},\dots, x_{n}\}\subset \phi^{c}\}=\prod_{i=1}^{n-2}
\sin^{2}\theta_{i}\times  \parallel \tilde{W}_{2}\parallel^{2} $
\item[(c)]
$P\{x_{1}\in \phi \vert \{x_{2},\dots,x_{n-1}\}\subset \phi^{c}\}=\parallel V_{1} \parallel^{2} + \parallel W_{1}\parallel^{2}$
\item[(d)]
\begin{equation}\label{pr2}
\begin{split}
& P\{x_{1}\in \phi \vert \{x_{2},\dots,x_{n}\}\subset \phi^{c}\}
= \\
& P\{x_{1}\in \phi, x_{n} \in \phi^{c}\vert \{x_{2},\dots,x_{n-1}\}\subset \phi^{c}\}\times \frac{P\{\{x_{2},\dots, x_{n-1}\}\subset \phi^{c}\}}{P\{\{x_{2},\dots, x_{n}\}\subset \phi^{c}\}}
\\
& =\big[  P\{x_{1}\in \phi\vert \{x_{2},\dots,x_{n-1}\}\subset \phi^{c}\}
\\
& \qquad
 - P\{x_{1}\in \phi, x_{n} \in \phi\vert \{x_{2},\dots,x_{n-1}\}\subset \phi^{c}\}\big] \\
& \quad \times \frac{P\{\{x_{2},\dots, x_{n-1}\}\subset \phi^{c}\}}{P\{\{x_{2},\dots, x_{n}\}\subset \phi^{c}\}}
\\
& =
\big[\parallel V_{1}\parallel^{2} + \parallel W_{1}\parallel^{2}
 -\parallel (V_{1},W_{1})\wedge (V_{2},W_{2})\parallel^{2}\big]
\times \frac{1}{\parallel\tilde{W}_{2}\parallel^{2}}.
\end{split}
\end{equation}
\end{enumerate}
From (a)-(d) an elementary computation 
gives the right side of (\ref{pp}) (note that  
$\parallel V_{2} \parallel^{2} + \parallel W_{2} \parallel^{2}+\parallel \tilde{W}_{2}\parallel^{2} = 1$). Indeed we get
\begin{equation}\label{pr3}
\begin{split}
& P\{\{x_{2},\dots, x_{n}\}\subset \phi^{c}\} \times 
P\{\{x_{2},\dots, x_{n-1}\}\subset \phi^{c}\}\\
& \times
\left[ P\{x_{1}\in \phi \vert \{x_{2},\dots,x_{n}\}\subset \phi^{c}\}
-P\{x_{1}\in \phi \vert \{x_{2},\dots,x_{n-1}\}\subset \phi^{c}\}\right]\\
& = \prod_{i=2}^{n-2}\sin^{4}\theta_{i}\left[
(\parallel V_{1} \parallel^{2} + \parallel W_{1} \parallel^{2})
(1 - \parallel \tilde{W}_{2}\parallel^{2})-
\parallel(V_{1},W_{1})\wedge (V_{2},W_{2})\parallel^{2}\right]
\\
& =  \prod_{i=1}^{n-2}\sin^{4}\theta_{i}\vert \langle(V_{1},W_{1}),(V_{2},W_{2})\rangle\vert^{2}.
\end{split}
\end{equation}
To compute the left side of (\ref{pp}) denote 
$z_{1}^{0}=(\sin\theta_{1}V^{1}_{1},   \dots,  \sin \theta_{n-2}
V^{n-2}_{1})$,\\
 $z_{n}^{0}=(\sin\theta_{1}V^{1}_{2},   \dots,  \sin \theta_{n-2}
V^{n-2}_{2})$
and
$\tilde{z}^{i}=(\cos \theta_{i}u^{i}_{1},  \dots,  \cos \theta_{i}u^{i}_{n-2},0)^{t}$.
Observe that
\begin{equation}\label{pr4}
 \langle z_{1},z_{n}  \rangle +
\sum_{k=1}^{n-2}(-1)^{k}\sum_{2\leq i_{1}<\dots < i_{k}\leq n-1}
  \langle z_{1}\wedge(\bigwedge_{j=1}^{k} z_{i_{j})},
 z_{n}\wedge(\bigwedge_{j=1}^{k} z_{i_{j}})  \rangle =A + B
 \end{equation}
 with
\begin{equation}\label{pr4'}
A = \langle z_{1}^{0},z_{n}^{0}  \rangle +
\sum_{k=1}^{n-2}(-1)^{k}\sum_{2\leq i_{1}<\dots < i_{k}\leq n-1}
  \langle z_{1}^{0}\wedge(\bigwedge_{j=1}^{k} z_{i_{j}}),
 z_{n}^{0}\wedge(\bigwedge_{j=1}^{k} z_{i_{j}})  \rangle
  \end{equation}
 and
 $$
 B=<W_{1},W_{2}>\big(1 +
\sum_{k=1}^{n-2}(-1)^{k}\sum_{1\leq i_{1}<\dots < i_{k}\leq n-2}
\parallel \bigwedge_{j=1}^{k}\tilde{z}^{i_{j}} \parallel^{2}\big)
$$
Obviously 
$\parallel \bigwedge_{j=1}^{k}\tilde{z}^{i_{j}} \parallel^{2}= 
\prod_{j=1}^{k}\cos^{2}\theta_{i_{j}},
$
thus
$$
1 +
\sum_{k=1}^{n-2}(-1)^{k}\sum_{1\leq i_{1}<\dots < i_{k}\leq n-2}
\parallel \bigwedge_{j=1}^{k}\tilde{z}^{i_{j}} \parallel^{2}
=\prod_{i=1}^{n-2}(1-\cos^{2}\theta_{i})=\prod_{i=1}^{n-2
}\sin^{2}\theta_{i}
$$
and consequently 
\begin{equation}\label{pr5}
B=<W_{1},W_{2}>\prod_{i=1}^{n-2}\sin^{2}\theta_{i}.
\end{equation}\\
In order to compute A introduce  the notations \\
$\tilde{z}^{i,l}=(\cos \theta_{i}u^{i}_{1},  \dots,  \cos \theta_{i}u^{i}_{n-2},\sin\theta_{i}V^{i}_{l})^{t}$, $l=1,2$ and $i=1,\dots,n-2$.\\
A little tought provides that for $k\geq1$
\begin{equation}\label{pr6}
\begin{split}
& \sum_{2\leq i_{1}<\dots < i_{k}\leq n-1}
  \langle z_{1}^{0}\wedge(\bigwedge_{j=1}^{k} z_{i_{j})},
 z_{n}^{0}\wedge(\bigwedge_{j=1}^{k} z_{i_{j}})  \rangle
 \\
 &=
\sum_{1\leq i_{1}<\dots < i_{k+1}\leq n-2}
 \bigg[  \langle \bigwedge_{j=1}^{k+1} \tilde{z}^{i_{j},1},
 \bigwedge_{j=1}^{k+1} \tilde{z}^{i_{j},2}   \rangle 
 - \parallel \bigwedge_{j=1}^{k+1} \tilde{z}^{i_{j}}\parallel^{2}\bigg]. 
 \end{split}
\end{equation}
Moreover we have 
\begin{equation}\label{pr7}
\begin{split}
&\bigwedge_{j=1}^{k} \tilde{z}^{i_{j},l}
=  \bigwedge_{j=1}^{k} \big( \tilde{z}^{i_{j}}+(e(n-2),\sin\theta_{i_{j}}V^{i_{j}}_{l})^{t}\big)\\
& =\bigwedge_{j=1}^{k} \tilde{z}^{i_{j}} + \sum_{j=1}^{k}(-1)^{j+1}(e(n-2),\sin\theta_{i_{j}}V^{i_{j}}_{l})^{t})\wedge(\bigwedge_{s=1, s\neq j}^{k} \tilde{z}^{i_{s}} ).
 \end{split}
\end{equation}
From orthogonality properties of relevant multivectors we obtain from (\ref{pr7})
\begin{equation}\label{pr8}
\begin{split}
&  \langle \bigwedge_{j=1}^{k} \tilde{z}^{i_{j},1},
 \bigwedge_{j=1}^{k} \tilde{z}^{i_{j},2}   \rangle 
 - \parallel \bigwedge_{j=1}^{k} \tilde{z}^{i_{j}}\parallel^{2}\\
& =  \sum_{j=1}^{k}\langle
 (e(n-2),\sin\theta_{i_{j}}V^{i_{j}}_{1})^{t}\wedge(\bigwedge_{s=1, s\neq j}^{k} \tilde{z}^{i_{s}} ),
 \\
 & \qquad \qquad
 (e(n-2),\sin\theta_{i_{j}}V^{i_{j}}_{2})^{t}\wedge(\bigwedge_{s=1, s\neq j}^{k} \tilde{z}^{i_{s}} )
  \rangle \\
  & = \sum_{j=1}^{k}V^{i_{j}}_{1}V^{i_{j}}_{2}\sin^{2}\theta_{i_{j}}
  \parallel \bigwedge_{s=1, s\neq j}^{k} \tilde{z}^{i_{s}}\parallel^{2}\\
  & =\sum_{j=1}^{k}V^{i_{j}}_{1}V^{i_{j}}_{2}\sin^{2}\theta_{i_{j}}
  \prod_{s=1, s\neq j}^{k}\cos^{2}\theta{i_{s}}.
 \end{split}
\end{equation}
From (\ref{pr4'}), (\ref{pr6}) and (\ref{pr8}) an elementary computation gives 
\begin{equation}\label{pr9}
A=\sum_{i=1}^{n-2}V^{i}_{1}V^{i}_{2}\sin^{2}\theta_{i}
\prod_{j=1,j\neq i}^{n-2}(1- \cos^{2}\theta_{j})
=
\langle V_{1},V_{2}\rangle\prod_{i=1}^{n-2
}\sin^{2}\theta_{i}
\end{equation}
and with (\ref{pr5}) 
\begin{equation}\label{pr10}
A+B=\langle (V_{1},W_{1}),(V_{2},W_{2})\rangle\prod_{i=1}^{n-2
}\sin^{2}\theta_{i}.
\end{equation}
This and  (\ref{pr3}) prove Proposition \ref{prop2}. Note that
from (\ref{pr10}) we get that (\ref{pr4}) is identified as a scalar product. 
\end{proof}
For further results by using the CSD as well as for some extensions of   Proposition \ref{prop2} see \cite{gol}.
\section{The BK inequality for increasing events generated by disjoint sets}
Let $\mathfrak{A}$, $\mathfrak{B}$ $\subset 2^{\mathcal{N}}$,
$\mathcal{N}=\{1,\dots,N\}$, be a pair of increasing events and suppose (obviously) that 
 $\emptyset \notin \mathfrak{A}\cup \mathfrak{B}$.
The events being increasing, there exist

 two minimal sets 
$ S_{1}= S(\mathfrak{A})=\{A_{i},i=1,\dots, n_{1}\}\subset \mathfrak{A}$ and
$ S_{2}=S(\mathfrak{B})=\{B_{i},i=1,\dots, n_{2}\}\subset \mathfrak{B} $
such that
\begin{enumerate}
\item  $ A\in \mathfrak{A} \quad \Leftrightarrow \exists A_{i} \quad such \quad that \quad A\supset A_{i} $,
\item $ B\in \mathfrak{B} \quad \Leftrightarrow \exists B_{i} \quad such \quad that \quad B\supset B_{i}. $
\end{enumerate}
The sets $A_{i}$,  $B_{i}$ are minimal in the sense that none of 
$A\in \mathfrak{A} $ (resp. $B\in \mathfrak{B} $) is strictly included
in $A_{i}$ (resp. in $B_{i}$).\\ 
Consider now a basic determinantal process $\phi$ on $\mathcal{N}$. In the particular case when
\begin{equation}
A \cap B =\emptyset \quad for \quad  all \quad  A\in S _{1} \quad and \quad \quad B\in S_{2}
\end{equation}
 we have at once $P\{\phi \in \mathfrak{A}\cap\mathfrak{B}\}
= P\{\phi \in \mathfrak{A}\circ\mathfrak{B}\}$ ans thus the 
BK-inequality (\ref{BK}) becomes
\begin{equation}\label{NA} 
P\{\phi \in \mathfrak{A}\cap\mathfrak{B}\}
\leq 
P\{\phi \in \mathfrak{A}\}
\times
P\{\phi \in \mathfrak{B}\}
\end{equation}
which is a negative association inequality. R.Lyons proved in \cite{lyo}, \cite{lyon} 
that determinantal processes have negative association, meaning that (\ref{NA}) 
is fulfilled. 
\\
In the general situation it is helpful 
 to reformulate BK inequality (\ref{BK}) as follows. 
\begin{prop}\label{prBK}
The inequality (\ref{BK}) is satisfied if and only if
\begin{equation}\label{p2}
\begin{split}
P\{\phi \notin \mathfrak{A}\cup\mathfrak{B}\}
\leq 
P\{\phi \notin \mathfrak{A}\}
\times
P\{\phi \notin \mathfrak{B}\} & + 
P\{\phi \in \mathfrak{A}\cap\mathfrak{B}\}\\
& - P\{\phi \in \mathfrak{A}\circ\mathfrak{B}\}.
\end{split}
\end{equation}
\end{prop}
\begin{proof}
Observe that
\begin{equation}
\nonumber
\begin{split}
P\{\phi \notin \mathfrak{A}\cup\mathfrak{B}\}
 & =1 - P\{\phi \in \mathfrak{A}\cup\mathfrak{B}\}
=1-  P\{\phi \in \mathfrak{A}\} -  P\{\phi \in \mathfrak{B}\}\\
& \qquad +P\{\phi \in \mathfrak{A}\cap\mathfrak{B}\}\\
& =P\{\phi \notin \mathfrak{A}\}
\times
P\{\phi \notin \mathfrak{B}\} - P\{\phi \in \mathfrak{A}\}
\times
P\{\phi \in \mathfrak{B}\} \\
& \qquad + P\{\phi \in \mathfrak{A}\cap\mathfrak{B}\}. 
\end{split}
\end{equation}
Thus 
\begin{equation}
\nonumber
P\{\phi \notin \mathfrak{A}\cup\mathfrak{B}\} -
P\{\phi \notin \mathfrak{A}\}
\times
P\{\phi \notin \mathfrak{B}\} -P\{\phi \in \mathfrak{A}\cap\mathfrak{B}\} + P\{\phi \in \mathfrak{A}\circ\mathfrak{B}\} \leq 0
\end{equation}
if and only if 
\begin{equation}
\nonumber
P\{\phi \in \mathfrak{A}\circ\mathfrak{B}\} - P\{\phi \in \mathfrak{A}\}
\times
P\{\phi \in \mathfrak{B}\} \leq 0.
\end{equation}
\end{proof}
Suppose now that $\mathfrak{A}=\mathfrak{B}$. Formula 
(\ref{p2}) becomes
 \begin{equation}\label{p4}
 P\{\phi \notin \mathfrak{A}\}
\leq 
P\{\phi \notin \mathfrak{A}\}^{2}
+ 
P\{\phi \in \mathfrak{A}\}
- P\{\phi \in \mathfrak{A}\circ\mathfrak{A}\}.
\end{equation}

If the sets of $S(\mathfrak{A})=\{A_{1},\dots,A_{n}\}$ are disjoint, that is if
  
 $A_{i}\cap A_{j}=\emptyset  $ for all $ i\neq j $,

then
$$ \{\phi \in \mathfrak{A} 
\backslash(\mathfrak{A}\circ\mathfrak{A})\}=
\bigcup_{i=1}^{n}\lbrace A_{i}\subset \phi, 
A_{j}\not \subset\phi, \forall j \neq i\rbrace.
 $$
Therefore
 
\begin{eqnarray*}\label{p5}
P\{\phi \in \mathfrak{A}\}
&-& P\{\phi \in \mathfrak{A}\circ\mathfrak{A}\}
= 
P\{\mathfrak{A}\backslash(\mathfrak{A}\circ\mathfrak{A})\}\\
&=& \sum_{i=1}^{n}P\lbrace A_{i}\subset \phi, 
A_{j}\not \subset \phi, \forall j \neq i\rbrace \\
 &=& \sum_{i=1}^{n}\big[ P\lbrace 
A_{j}\not \subset \phi, \forall j \neq i\rbrace
- P\lbrace A_{i}\not \subset\phi, \forall i=1,\dots,n\rbrace 
\big]\\
&=& \sum_{i=1}^{n}P\lbrace 
A_{j}\not \subset \phi, \forall j \neq i\rbrace
- nP\lbrace \phi \notin \mathfrak{A}\rbrace
\end{eqnarray*}
 
and formula  (\ref{p4}) takes the following form
\begin{eqnarray}\label{p6}
(n+1)P\{\phi \notin \mathfrak{A}\}
\leq
P\{\phi \notin \mathfrak{A}\}^{2}
 +   
 \sum_{i=1}^{n}P\lbrace 
A_{j}\not \subset \phi, \forall j \neq i\rbrace. 
\end{eqnarray}
Fix  now $n_{0}\geq 2$ and suppose that the Conjecture \ref{Con} is fulfilled for all $2\leq n\leq n_{0} $.
\begin{lemma}\label{lp1}
Under this hypothesis, for all  
$A_{i}$, $i=1,\dots,n$,  disjoint subsets of $\{1,\dots, N\}$ with 
$2\leq n\leq n_{0} $, such that $P\lbrace  A_{i} \not \subset \phi, \forall i=1,\dots,n\rbrace>0$
we have
\begin{equation}\label{ep1}
P\lbrace  A_{i} \not \subset \phi, \forall i=1,\dots,n\rbrace^{n-1}
\leq 
\prod_{i=1}^{n}P\lbrace  A_{j} \not \subset \phi, \forall j\neq i\rbrace.
\end{equation}
\end{lemma}
\begin{proof} 
For $n=2$ the inequality (\ref{ep1}) is the well-known  correlation inequality. For $n > 2$ applying (\ref{C}) we get
\begin{equation}\label{ep2} 
\begin{split}
& \prod _{k=2}^{n}P\lbrace 
A_{k}\not \subset \phi \mid  A_{j} \not \subset \phi,
 \forall j \neq k \rbrace 
 \leq 
\prod_{k=2}^{n}P\lbrace 
A_{k} \not \subset \phi \mid  A_{j} \not \subset \phi,
 \forall j \neq 1, k \rbrace\\
 & \Longleftrightarrow \\
 &
\dfrac{P\lbrace 
A_{i}\not \subset \phi , \forall i=1,\dots,n \rbrace ^{n-1}}{\prod_{i=1}^{n}P\lbrace 
A_{j} \not \subset \phi, \forall j \neq i\rbrace} 
\leq 
\dfrac{P\lbrace 
A_{i} \not \subset \phi, \forall i \neq 1 \rbrace ^{n-2}}{\prod_{k=2}^{n}P\lbrace 
A_{j} \not \subset \phi, \forall j \neq 1,k\rbrace} 
 \end{split}.
\end{equation}
and thus Lemma \ref{lp1} follows by induction.
\end{proof}
We will need the following elementary lemma. Its proof being trivial we omit it.

\begin{lemma}\label{lp2}
For all $ 0<a\leq 1 $ and $n>0$ we have 
\begin{equation}
(n+1) -a - na^{-\frac{1}{n}}\leq 0.
\end{equation}
\end{lemma}

\begin{theorem}\label{t1g}
Let  $\mathfrak{A}$ be an increasing event generated by disjoint sets
$ A_{1},\dots,A_{n}$. Suppose that the Conjecture \ref{Con} holds.   We have
\begin{equation}\label{B7g} 
P\{\phi \in \mathfrak{A}\circ\mathfrak{A}\}
\leq 
P\{\phi \in \mathfrak{A}\}^{2}.
\end{equation}
\end{theorem}
\begin{proof} We have to prove (\ref{p6}). By Lemma \ref{lp2}
 we obtain
\begin{eqnarray}\label{B8g}
(n+1)P\{\phi \notin \mathfrak{A}\}
\leq P\{\phi \notin \mathfrak{A}\}^{2}
+ nP\{\phi \notin \mathfrak{A}\}^{\frac{n-1}{n}}.
\end{eqnarray}
Lemma \ref{lp1} implies that 
\begin{equation}
P\{\phi \notin \mathfrak{A}\}^{n-1}
=P\lbrace  A_{i} \not \subset \phi, \forall i=1,\dots,n\rbrace^{n-1}
\leq \prod_{i=1}^{n}P\lbrace  A_{j} \not \subset \phi, \forall j\neq i\rbrace
\end{equation}\label{B9g}
so it remains to apply the 
geometric-arithmetic mean inequality
\begin{equation}\label{B10g}
n\prod_{i=1}^{n}P\lbrace  A_{j} \not \subset \phi, \forall j\neq i\rbrace^{\frac{1}{n}}
\leq \sum_{i=1}^{n}P\lbrace  A_{j} \not \subset \phi, \forall j\neq i\rbrace
\end{equation}
to obtain (\ref{p6}) as desired.
\end{proof}
\begin{remark}
Consider an event $\tilde{S}=\{ D_{1},\dots,D_{n_{0}}\}\subset  2^{\mathcal{N}}$ of disjoint sets such that 
$P\{D \not \subset \phi, \forall D \in \tilde{S} \}>0$.
Denote $\psi = \{\phi \vert D \not \subset  \phi,  \forall D \in \tilde{S}\}$. If  the Conjecture \ref{Con} holds then it is obvious that the inequality
 (\ref{C}) is also satisfied for the conditioned process $\psi$
 provided that the sets occuring in (\ref{C}) are disjoint from those  in $\tilde{S}$. Consequently, if  $\mathfrak{A}$ is an increasing event generated by disjoint sets
$ A_{1},\dots,A_{n} $ such that $ A_{i}\cap D = \emptyset$
for all $i=1,\dots,n$ and $D\in \tilde{S}$ then we obtain
\begin{equation}\label{B7g+} 
P\{\psi \in \mathfrak{A}\circ\mathfrak{A}\}
\leq 
P\{\psi \in \mathfrak{A}\}^{2}.
\end{equation}
\end{remark}
Let   $S_{1}=\{A_{i},i=1,\dots, n_{1}\}$, $S_{2}=\{B_{i},i=1,\dots, n_{2}\}$ and $S=\{C_{i},i=1,\dots, n_{3}\}$ be events
such that all sets in  $S_{1}\cup S_{2}\cup S \subset 
2^{\mathcal{N}}$ are pairwise disjoint.  
\begin{theorem}\label{t2g}
Suppose that the Conjecture \ref{Con} holds. Then for increasing events $\mathfrak{A}$, $\mathfrak{B}$ such that $S(\mathfrak{A})=S_{1}\cup S$ and $S(\mathfrak{B})=S_{2}\cup S$ we have
\begin{equation}\label{B19g} 
P\{\psi \in \mathfrak{A}\circ\mathfrak{B}\}
\leq 
P\{\psi \in \mathfrak{A}\}\times P\{\psi \in \mathfrak{B}\}. 
\end{equation}
where $\psi =\{ \phi \vert D \not \subset \phi,  \forall D \in \tilde{S}\}$ and all sets in $S_{1}\cup S_{2}\cup S \cup \tilde{S}\subset 
2^{\mathcal{N}}$ are pairwise disjoint.  
\end{theorem}
\begin{proof}
The proof proceeds by induction using Theorem \ref{t1g} and the following lemma.
\begin{lemma}\label{lp3}
Fix $S_{1}$, $S_{2}$, $S$ and  and suppose that BK inequality (\ref{B19g}) is fulfilled for all conditioned processes $\psi$ subjected to the conditions of Theorem \ref{t2g}. Fix $A \subset \mathcal{N}$, $A\neq \emptyset$ such that $A\cap A' = \emptyset$ for all $A' \in S_{1}\cup S_{2}\cup S $. Denote by 
$  \mathfrak{\tilde A}=\sigma \{A,\mathfrak{A}\}$ the increasing event
generated by $A
$ and $ \mathfrak{A}$. Then, the
BK inequality
  \begin{equation}\label{B21g}
P\{\psi \in \mathfrak{\tilde A}\circ\mathfrak{B}\}
\leq 
P\{\psi \in \mathfrak{\tilde A}\}
\times
P\{\psi \in \mathfrak{B}\}
\end{equation}
is satisfied for all conditioned processes 
 $\psi = \{\phi \vert D \not \subset \phi \quad \forall D \in \tilde{S}\}$ such that  all sets of $S(\mathfrak{\tilde A})\cup S_{2}\cup S \cup \tilde{S}\subset 
2^{\mathcal{N}}$ are pairwise disjoint.  
\end{lemma}
\begin{proof} 
By  (\ref{p2}) we may suppose that $$P\{A' \not \subset \psi, \forall A' \in  S(\mathfrak{\tilde A})\cup S_{2}\cup S \cup \tilde{S}\}>0.$$
We have
\begin{equation}\label{B22g}
\lbrace \psi \in \mathfrak{A}\cap\mathfrak{B}
\setminus \mathfrak{A}\circ\mathfrak{B}\rbrace=
\cup_{C\in S}\lbrace C \subset \psi, 
A' \not \subset \psi, \forall A' \in S_{1}\cup S_{2}\cup S
\quad A'\neq C\rbrace 
\end{equation}
and
\begin{equation}\label{B23g}
\lbrace \psi \in\mathfrak{\tilde A}\cap\mathfrak{B}
\setminus \mathfrak{\tilde A}\circ\mathfrak{B}\rbrace=
\cup_{C\in S}\lbrace C \subset \psi, A\not \subset \psi,  
A' \not \subset \psi, \forall A' \in S_{1}\cup S_{2}\cup S
\quad A'\neq C\rbrace. 
\end{equation}
Formulas (\ref{p2}) and (\ref{B23g}) imply that the BK inequality (\ref{B21g}) can be written as follows
\begin{equation}\label{B24g}
\begin{split}
P\{ A\not \subset \psi, & A' \not \subset \psi, \forall A' \in S_{1}\cup S_{2}\cup S\}
\leq 
P\{A\not \subset \psi, A' \not \subset \psi, \forall A' \in S_{1}\cup S\}\\
& \qquad \qquad
\times P\{ A' \not \subset \psi, \forall A' \in S_{2}\cup S\} \\
& + 
\sum_{C\in S}P\lbrace C \subset \psi, A\not \subset \psi, 
A' \not \subset \psi, \forall A' \in S_{1}\cup S_{2}\cup S
\quad A'\neq C \rbrace 
\end{split}
\end{equation}
or, introducing the process $ \psi_{0} = \{\phi \mid  A \not\subset \phi,  A' \not\subset \phi \quad \forall A' \in \tilde{S}\}$, as 
\begin{equation}\label{B25g}
\begin{split}
P\{ & A' \not \subset \psi_{0}, \forall A' \in S_{1}\cup S_{2}\cup S\}
\leq 
P\{A' \not \subset \psi_{0}, \forall A' \in S_{1}\cup S\}\\
& \qquad \qquad
\times P\{ A' \not \subset \psi, \forall A' \in S_{2}\cup S\} \\
& + 
\sum_{C\in S}P\lbrace C \subset \psi_{0}, 
A' \not \subset \psi, \forall A' \in S_{1}\cup S_{2}\cup S
\quad A'\neq C \rbrace.
\end{split}
\end{equation}
The stated hypotheses imply that
\begin{equation}\label{B26g}
\begin{split}
P\{ & A' \not \subset \psi_{0}, \forall A' \in S_{1}\cup S_{2}\cup S\}
\leq 
P\{A' \not \subset \psi_{0}, \forall A' \in S_{1}\cup S\}\\
& \qquad \qquad
\times P\{ A' \not \subset \psi_{0}, \forall A' \in S_{2}\cup S\} \\
& + 
\sum_{C\in S}P\lbrace C \subset \psi_{0}, 
A' \not \subset \psi, \forall A' \in S_{1}\cup S_{2}\cup S
\quad A'\neq C \rbrace. 
\end{split}
\end{equation}
It is easy to see that Conjecture \ref{Con} implies the inequality
\begin{equation}\label{B27g} 
P\{ A' \not \subset \psi_{0}, \forall A' \in S_{2}\cup S\}
\leq
P\{ A' \not \subset \psi, \forall A' \in S_{2}\cup S\}
\end{equation}
and thus by (\ref{B26g}) and (\ref{B27g}) we obtain (\ref{B25g}) which finish the proof of Lemma \ref{lp3}.
\end{proof}
Starting from  (\ref{B7g+})  
and  applying step by step the Lemma \ref{lp3}  Theorem \ref{t2g} follows. 
\end{proof}
\section{The BK inequality for increasing events $\mathfrak{A}$, $\mathfrak{B}$ generated by simple points}
As mentioned in the Introduction the inequality (\ref{C}) is satisfied when  the occuring sets are reduced to being simple points. This follows easily, for example, from Proposition \ref{prop1}. Therefore  Theorem \ref{t2g} implies
\begin{theorem}\label{t2}
Let  $\mathfrak{A}$, $\mathfrak{B}$ be increasing events generated
by simple points. The BK inequality
\begin{equation}\label{B19} 
P\{\phi \in \mathfrak{A}\circ\mathfrak{B}\}
\leq 
P\{\phi \in \mathfrak{A}\}\times P\{\phi \in \mathfrak{B}\} 
\end{equation}
is then satisfied for all determinantal discrete processes $ \phi $ associated to  sets of orthonormal vectors of $ \mathbb{C}^{N} $. 
\end{theorem}
\begin{remark}\label{R4}
For sets reduced to being  simple points the key inequality
(\ref{ep1}) can be seen from the point of view given by the CSD. Indeed, 
consider the CSD in the case I applied to   
 $J=\{x_{1},\dots,x_{n}\}$ and, according, let  
$ v^{j}=(v_{1}^{j},\dots,v_{n}^{j})^{t}$,
$v_{i}^{j}=(\sin\theta_{j})u_{i}^{j}$, $i,j=1,\dots,n$ 
be
the vectors  such that
\begin{equation}
P\{\{x_{i_{1}},\dots,x_{i_{k}}\} \subset \phi^{c}\}=
\parallel \bigwedge_{j=1}^{k}v_{i_{j}}\parallel ^{2}
\end{equation}
for all $\{x_{i_{1}},\dots,x_{i_{k}}\} \subset J$.
Denote  
\begin{equation}\label{B1}
\tilde{v}_{i}=\bigwedge\limits_{j\neq i}v_{j}=(\tilde{v}_{i}^{1},\dots,\tilde{v}_{i}^{n}) \in \mathbb{C}^{n}
, i=1,\dots,n
\end{equation}
where \begin{equation}\label{B5}
\tilde{v}_{i}^{j}= \prod_{k\neq j}\sin\theta_{k}\times \tilde{u}_{i}^{j}
\end{equation}
and
$\tilde{u}_{i}^{j}$ is the $(i,n-j+1$)-minor of the unitary matrix $U=(u^{j}_{i})_{i,j=1,\dots,n}$. 
By (\ref{itt}) we obtain:
 
 \begin{equation}\label{B10}
\begin{split}
P\{x_{i}\in \phi^{c}, i=1,\dots,n\}^{n-1} & =\prod_{i=1}^{n}(\sin\theta_{i})^{2(n-1)} \\
& =
\parallel \bigwedge_{i=1}^{n}v_{i}\parallel ^{2(n-1)}
\\
&=\parallel \bigwedge_{i=1}^{n}\tilde{v}_{i}\parallel^{2}
\\
& \leq \prod_{i=1}^{n}\parallel \tilde{v}_{i}\parallel^{2}
=\prod_{i=1}^{n}P\lbrace 
x_{j}\in \phi^{c}, \forall j \neq i\rbrace .
\end{split} 
\end{equation}
\end{remark}
\begin{remark}\label{R5}
It was pointed out to us that for an increasing event $\mathfrak{A}$ generated by simple points $S= \{x_{1},\dots,x_{n}\}$ the inequality (\ref{B7g}) which can be read as 
\begin{equation}\label{B11}
P\{\vert S\cap \phi\vert\geq 2\}\leq P\{\vert S\cap \phi\vert\geq 1\}^{2}
\end{equation}
can be obtained also by a direct computation from formula (\ref{it1}) of Proposition \ref{prop1} and, moreover, if one consider the product measure $\mu=\otimes_{i=1}^{n}((\cos^{2}\theta_{i})\delta_{1} + (\sin^{2}\theta_{i})\delta_{0})$ on the product space $E=\{0,1\}^{n}$ and increasing events 
$$\mathfrak{A}_{i}=
\{a=(a_{j})\in E \quad such \quad that \quad \sum_{j=1}^{n}a_{j} \geq i\},$$
$i=0,\dots,n$, then the formulas (\ref{it1}) imply that $P\{\vert S\cap \phi\vert\geq i\}
= \mu (\mathfrak{A}_{i})$. From the original J. van den Berg and  H.Kesten Theorem (3.3) of \cite{BK} we get that
\begin{equation}\label{B12}
P\{\vert S\cap \phi\vert\geq i+j\}\leq P\{\vert S\cap \phi\vert\geq i\}\times P\{\vert S\cap \phi\vert\geq j\}, \quad 2\leq i+ j \leq n.
\end{equation} 
Furthermore, note that by Remark \ref{R2} 
the inequalities (\ref{B11}) and (\ref{B12}) are still valid for general determinantal processes (both discrete and continuous) taking for S a Borel set.

\end{remark}

\section{Extensions and concluding remarks}
Theorem \ref{t2} can be easily extended in the setting of general discrete determinantal processes. From the construction  given in paragraphs 2.2. of \cite{lyon}, which start from the basic processes, it follows at once that Theorems 1-2 are valid 
(the generated sets $S(\mathfrak{A})$ and $S(\mathfrak{B})$ being finite
or infinite)
for determinantal point processes defined on denumerables sets 
$\mathcal{E}$ and associated to closed subspaces of $l^{2}(\mathcal{E})$. Now, let $\phi$ be a such process on $\mathcal{E}$. Fix 
$\mathcal{F} \subset \mathcal{E}$ and consider the proces 
\begin{equation}\label{Q} 
\psi= \phi \cap \mathcal{F}.
\end{equation} 
Let $\mathfrak{A}$, 
$\mathfrak{B}\subset 2^{\mathcal{F}}$,  $\tilde{\mathfrak{A}}$, 
$\tilde{\mathfrak{B}}\subset 2^{\mathcal{E}}$ be the increasing events generated respectively by $S_{1}=S(\mathfrak{A})=S(\tilde{\mathfrak{A}}) \subset \mathcal{F}$
and $S_{2}=S(\mathfrak{B})=S(\tilde{\mathfrak{B}}) \subset \mathcal{F}$. The BK inequalities for $\phi $, $\tilde{\mathfrak{A}}$,  $\tilde{\mathfrak{B}}$
and $\psi$, $\mathfrak{A}$, 
$\mathfrak{B}$,  involve only generating sets $S_{1}$, $S_{2}$. Consequently Theorem \ref{t2} is valid for $\psi$ as well. To finish just note that, by paragraph 2.2. of \cite{lyon},  discrete determinantal processes  associated to  positive contractions (the general case) are of the form (\ref{Q}).\\
By the transference principle (\cite{lyon} 3.6.) Theorem \ref{t2} could also be extended to the continuous case but this is of little use due the fact that in the continuous setting the intensity measures related to determinantal processes of interest are of diffusive type  which implies that $P\{x \in \phi \}=0$ for  points $x$ (however, as mentioned in  Remark \ref{R5}, inequalities (\ref{B11}) and (\ref{B12}) still hold).
 
\section*{Acknowledgement}
I thank the anonymous referees for their constructive comments, especially for a pertinent question about the validity of
the Conjecture \ref{Con} for the non-determinantal point processes, which led to look at the process  described in Remark \ref{R0}.

\end{document}